\input amstex
\documentstyle{amsppt}
\magnification=\magstep1
\vcorrection{-8mm}
\NoRunningHeads
\define\bs{\bigskip}
\define\ms{\medskip}
\define\hook{\hookrightarrow}
\redefine\F{\Cal F}
\topmatter
\title GENERAL THEOREM ON INTERPOLATION OF COMPACT OPERATORS
\endtitle
\author Evgeniy Pustylnik \\\\ Technion, Haifa 32000, Israel \\
{\it E-mail: evg\@gmail.com}
\endauthor
\keywords Compact operators, interpolation, basic sequences
\endkeywords
\subjclass 46B70 \endsubjclass
\abstract
We prove an abstract theorem on keeping the compactness property of a linear operator after interpolation in Banach spaces. No analytical presentation of operators, spaces and interpolation functor is required. We use only some little-known properties of compact sets and various facts about bases and basic sequences (with detailed references to the monograph ``Bases in Banach spaces", \!Vol.\,I--II  by I.M.Singer). Therefore the results are applicable to arbitrary spaces and any interpolation functor, {\sl including the complex method}. The ``two-sided" compactness is also mentioned at the end of this paper as a mere corollary.
\endabstract
\endtopmatter

\document

{\bf 1.} In this paper we show that the interpolation of ``one-sided" compactness of arbitrary linear operators
is possible for any Banach couples and every interpolation functor. We start our investigations with
embedding operators, since the consequent pass to arbitrary operators was already considered
by M. Cwikel, N. Krugljak and M. Masty\l o in [1]\ (see {\it Appendix}\,\ at the end of this paper).

So, let $(A_0,A_1),\,(B_0,B_1)$ be two nontrivial Banach couples with an interpolation functor $\F$,
defined on them. Let $A_0\hook B_0,\,A_1\hook B_1$ and $J$ be the corresponding embedding operator. Thus,
due to interpolation, $J$ gives an embedding $\F(A_0,A_1)=A\hook B=\F(B_0,B_1)$. The problem is to prove
compactness of the last embedding if such is the embedding $J: A_1\to B_1$ (one-sided compactness).

Of course, the one-sided compactness of embedding for an intermediate space may be impossible if this space is
not sufficiently distant from the second space of the couple, where the compactness is not given. This imposes
some additional restrictions on the interpolation functor, which must contain an implicit measure of admissible
distance. For example, those restrictions are shown as necessary for some cases of real interpolation in the
papers [2], [3] etc.

In our paper such a restriction will be used in the following form:\ {\sl there exists a function $W(\alpha,\beta),\
\alpha,\beta>0$, with $\lim_{\beta\to 0}W(\alpha,\beta)=0$ when $\sup\alpha<\infty$, such that
$$ \|T\|_{A\to B}\le W(\|T\|_{A_0\to B_0},\|T\|_{A_1\to B_1}) \tag 1 $$ for any linear operator $T:A_0+A_1\to
B_0+B_1$}. For example, in the case of a complex interpolation functor $\F_{[\theta]},\ 0<\theta<1$, such a
function exists in the form $W(\alpha,\beta)=\alpha^{1-\theta}\beta^\theta$. \ms

\proclaim{Main Theorem} Under conditions on the spaces $(A_0,A_1),\,(B_0,B_1)$ and the functor $\F$ posed above,
the embedding $A\hook B$ is compact. \endproclaim

We shall use the following standard way of the proof. We suppose, on the contrary, that the embedding $A\hook B$
is not compact. Then there exists a bounded infinitely dimensional sequence $(x_n)\subset A,\ x_n\neq0$, with no convergent
subsequence in $B$. It follows immediately from this that the sequence of the norms $\|x_n\|_B$ is bounded from
below too, namely, $\inf\|x_n\|_B>0$. This fact can be proved exactly as in Theorem 4.2 from [4] due to non-compactness
of embedding $A\hook B$, since every unbounded from below set necessarily contains a sequence converging to zero.
The infinite dimension of the sequence $(x_n)$ also follows from its non-compactness in $B$, because each set with
finite dimension is compact in any topology.

Together with the embedding inequality between norms in spaces $A$ and $B$, we obtain the following assertion.
\proclaim{Corollary} The assumption on non-compactness of embedding $A\hook B$ implies existence of an infinitely
dimensional bounded from above and below sequence $(x_n)$ with equivalent norms $\|x_n\|_A$ and $\|x_n\|_B$.
\endproclaim
Just one of such sequences will be the main object in our next discussions.\ms

{\bf 2.} Further on we shall permanently use a special ``reduction principle"\footnote"$^{(1)}$" {The term is new.}.
It is easy to see that our main problem will be solved if we prove it for some infinitely dimensional
subspaces of $A_0$ and $A_1$. Indeed, if the embedding of given spaces is compact, it remains such for
any smaller parts. Similarly, if the sequences of norms $\|x_n\|_A$ and $\|x_n\|_B$ are equivalent, this
equivalence is kept for corresponding norms of any subsequence $(x_{n_j})$. The ``reduction principle" is
intended for deleting all unnecessary parts of spaces (or sequences), because they
may prevent from use of some important arguments in the forthcoming proofs.

The main advantage of the ``reduction principle" in this paper is the possibility to deal with bases of all needed
subspaces. As known, a basis may not exist even in a separable Banach spaces. At the same time, as shown in [6]
(Theorem 1.1, page 48), any Banach space has a basic sequence. Recall that a sequence $(z_n)$ in a Banach space $E$ is
said to be basic one if $(z_n)$ is a Schauder basis of $[z_n]$, where $[z_n]$ means the closed linear subspace of $E$
spanned by the sequence $(z_n)$. (By the way, the books [5] and [6] by I. Singer will be our main reference for all
needed properties of the bases and basic sequences.) In conclusion: {\sl the ``reduction principle" allows us to replace
in future proofs any considered space $E$ by a suitable subspace $[z_n]$, obtaining a genuine basis instead of a basic
sequence}.

Yet more strong is an assertion from [6] (Corollary 1.5) about existence of a basic subsequence in arbitrary sequence,
having infinite dimension. In fact, this Corollary asserts only existence of a so-named ``block-basis sequence", but the
previous volume [5] already explained on the page 66 that the block-basis sequence is merely a special kind of basic sequences.
Now it must be clear why the infinite dimension of the sequence $(x_n)$ was especially stressed above, because only thus we
may state (after application of the ``reduction principle") that this sequence is a basis of the reduced space $A$.

Another useful possibility of the ``reduction principle" is transformation of any embedding $A\subset B$ into dense one,
reducing the space $B$ to the closure of $A$ in $B$. In our proof, this operation will be done for the space $\Delta A=
A_0\cap A_1$ which will be from now on regarded as dense in the space $A$. This allows us to replace the sequence $(x_n)$,
defined above, by a sequence $(y_n)$ from $\Delta A$ which is so close to $(x_n)$ as to have equivalent norms in $A$ and $B$
as well. Moreover, due to [5] (Theorem 10.1, page 93), the sequence $(y_n)$ remains to be basis in $A$. Hence a new sequence
completely replaces the former one, but, for convenience, we proceed to use the same notation $x_n$ instead of $y_n$, adding
only one new property that all $x_n$ belong to $\Delta A$.

By the definition of $\Delta A$, the sequence $(x_n)$ geometrically (as set) is the same in all three spaces $A,A_0,A_1$,
but the norms of elements may be different. However we may change these norms unless this involves any serious alteration
of the initial problem. For instance, we may take all norms $\|x_n\|_{A_0}=\|x_n\|_{A_1}=1$ and correspondingly obtain
(due to interpolation) some new values of $\|x_n\|_A\lesssim \|x_n\|_{\Delta A}=\max\{\|x_n\|_{A_0},\|x_n\|_{A_1}\}$
\footnote"$^{(2)}$" {As usual, the sign $\lesssim$ means an inequality with unknown coefficient like $a\le kb$.}. It
looks like multiplying all previous norms $\|x_n\|_A$ by some multipliers $\lambda_n$, bounded from above on the whole.
Moreover, the sequence $(\|x_n\|_A)$ proceeds to be bounded from below due to its non-compactness in the space $B$ as was
explained above while defining this sequence. Notice, in addition, that $A$, as any interpolation space, is intermediate in
the couple $(A_0,A_1)$, therefore the sequence $(x_n)$ cannot become unbounded in $A$ even obtaining new (bounded) norms
in $A_0,A_1$. Indeed, the linear independence of the elements $x_n$ with different indices implies that $\|x_n\|_{A_0+A_1}=
\min\{\|x_n\|_{A_0},\|x_n\|_{A_1}\}$ while $\|x_n\|_A\gtrsim\|x_n\|_{A_0+A_1}$.
Therefore the sequence of multipliers $(\lambda_n)$ is bounded both from above and from below and, by
Definition 3.2 and Theorem 3.2 from [5], pages 21@-22, $(x_n)$ remains a bounded basis of the space $A$ with the same topology.
By the way, the norms $\|x_n\|_B$ obtain the same multipliers as the norms $\|x_n\|_A$, and both norms remain equivalent. \ms

Let us mention two important details concerning the ``reduction principle". Excluding some elements $(x_n)$ from one of
spaces, we should simultaneously exclude them from two other spaces, having thus the same remaining sequence $(x_n)$ in any
of spaces at any time. The second remark is that we do not change the notations after any such actions as if we had all these
improved objects from the very beginning. This concerns the spaces $A_0,A_1,A$ and $B_0,B_1,B$ themselves as well. \ms

After reduction of the sequence $(x_n)$ in the space $A$ such that it became a basis, we begin processing it in the spaces
$A_0$ and $A_1$. Being a basis, this sequence is infinite dimensional and remains such in other spaces, because this property is
algebraical and does not depend on topology. Therefore we can reduce it till a genuine basis in the space $A_0$. At last, by the
next reduction, we could get an analogous basis in the space $A_1$. As shown in [5] (Proposition 4.1, page 26),
all previous bases remain be such even after reduction, so that we obtain a common basis $(x_n)$ in all three spaces.

The properties of the elements $x_n$ in the spaces $B,B_0,B_1$ can be obtained as consequences of the given
embeddings. Without loss of generality we may suppose that all embeddings are normalized, namely, that
$$ \|x\|_{B_0}\le\|x\|_{A_0}, \quad \|x\|_{B_1}\le\|x\|_{A_1}, \quad \|x\|_B\le\|x\|_A \tag 2 $$
for all admissible $x$.

\proclaim{Lemma} Let the embedding $A_1\hook B_1$ be compact and let $(x_n)$ be a basis of the space $A_1$ with unit
norms. Then $\lim_{n\to\infty}\|x_n\|_{B_1}=0$. \endproclaim

\demo{Proof} Suppose, on the contrary, that $\|x_n\|_{B_1}\not\rightarrow 0$. Then there exist infinitely many elements $x_{n_j}$
with norms in $B_1$ bigger than some positive number $\varepsilon_0$. Due to (2), these norms are less than $\|x_{n_j}\|_{A_1}=1$,
so that the sequence $(x_{n_j})$ has equivalent norms in the spaces $A_1$ and $B_1$. But, as shown in [4] (Lemma 4.1 and Theorem
4.2, page 324), this contradicts to compactness of the embedding $A_1\hook B_1$.

{\it A small remark:} Lemma 4.1 from [4] requires that the sequence $(x_{n_j})$ contains a basic subsequence in $B_1$.
But this sequence (geometrically, without comparing the norms) is a part of the basis in $A_1$ and thus has an infinite dimension.\qquad \qed \enddemo

{\bf 3.} For the next application of the ``reduction principle", we construct a special subsequence $(x_{n_j})$ of the sequence $(x_n)$ with
the property  $\|x_{n_{j+1}}\|_{B_1}\le\|x_{n_j}\|_{B_1}/2$ for all $j=1,2,\dots$ \ Due to Lemma, such a subsequence necessarily exists.
Making sufficient repetitions, we obtain
$$ \|x_{n_{j+k}}\|_{B_1}\le 2^{-k}\|x_{n_j}\|_{B_1}\quad \text{for any}\quad j,k\in\Bbb N. \tag 3 $$
Due to ``reduction principle", we may use now only the sequence $(x_{n_j})$, which keeps all needed properties of $(x_n)$: it is basic
in all three spaces $A, A_0, A_1$ (as subsequence of bases) and has equivalent norms in the spaces $A$ and $B$. This allows us to return
for this sequence to the previous notation $(x_n)$. The only novelty of the sequence $(x_{n_j})$, that could be lost, is the formula (3).
In former notations (without duple indices), it can be rewritten as
$$ \|x_{n+k}\|_{B_1}\le 2^{-k}\|x_n\|_{B_1}\quad \text{for any}\quad n,k\in\Bbb N. \tag 4 $$\ms

As a last component of proof, we use the connection between basis expansion and its partial sums ([5], Theorem 7.1, page 57).
Namely, if $(x_n)$ is a basis in a Banach space $A_1$ and $a_n$ are the numerical
coefficients, then there exists a constant $C$ such that, for any $N$ and any $(a_n)$, one has
$$ \|\sum_{n=1}^N a_nx_n\|_{A_1}\le C\|\sum_{n=1}^\infty a_nx_n\|_{A_1}. \tag 5 $$
In particular, it follows from this that
$$ \|a_kx_k\|_{A_1}=\|\sum_{n=1}^k a_nx_n-\sum_{n=1}^{k-1}a_nx_n\|_{A_1}\le 2C\|\sum_{n=1}^\infty a_nx_n\|_{A_1} \tag 6 $$ and
$$ \|\sum_{n=N}^\infty a_nx_n\|_{A_1}=\|\sum_{n=1}^\infty a_nx_n-\sum_{n=1}^{N-1} a_nx_n\|_{A_1}\le(C+1)\|\sum_{n=1}^\infty
   a_nx_n\|_{A_1}. \tag 7 $$
Let us notice that the same inequalities are hold after replacing $A_1$ by $A_0$. Moreover, since the basis $(x_n)$ is
common for both spaces, the constant $C$ depends here only on spaces $A_0, A_1$ themselves and can be chosen the same
(maximal of two) for either of them. \ms

{\bf 4.} Now, after all preliminaries, we can pass to the core of the proof. First of all, remark that both initial suppositions on
compactness of embedding $A_1\hook B_1$ and on equivalence of norms of the sequence $(x_n)$ in spaces $A$ and $B$ remained
true after all applications of the ``reduction principle". This allows us to retain the same notations of all spaces after
those reductions and to use for our proof all set of inequalities (2)--(7).

The inequality (1) needs a special discussion, because it relates to all spaces in their initial form, before all reductions.
In order to make a bridge between former and current states of all spaces and sequences, we define an operator $Q$, acting on
the initial space $A_0+A_1$ and equal to zero at all elements, deleted during all reductions. For the remaining elements, we set
$Qx=x$. In fact, $Q$ is a projection operator of $A_0+A_1$ onto the last subspace $[x_n]$, spanned by the remains of the sequence
$(x_n)$ after all reductions. Combining this operator with the embedding operator $J$, we obtain all corresponding reductions for
the space $B_0+B_1$. In the same manner we obtain that $TQ:A_0+A_1\to B_0+B_1$ for any operator $T$ which is defined only on reduced
couples. This allows us to use the inequality (1) for initial couples, replacing $T$ by $TQ$, or to use this inequality without
$Q$ for the same couples, obtained after reductions.

\demo{Proof of the Main Theorem} In the following proof all spaces are considered in the final form, after all reductions.
We also need a special kind of operator-remainders that are very useful in the theory of bases. Namely, if $x=\sum_{n=1}^
\infty a_nx_n$, then $R_nx:=\sum_{k=n}^\infty a_kx_k$. Next we recall the embedding operator $J:A_0+A_1\to B_0+B_1$ and
define the operators $P_n=R_nJ:A_0+A_1\to B_0+B_1,\ n=1,2,\dots$ Due to presence of the embedding operator, we obtain that
$P_n:A_0\to B_0$ and $A_1\to B_1$ separately, that is, these operators can be interpolated between couples $(A_0,A_1)$
and $(B_0,B_1)$.

Let us start with embedding $A_1\hook B_1$. Consider
$$ \|P_n\|_{A_1\to B_1}=\sup_{\|x\|_{A_1}\le1}\|P_nx\|_{B_1}. $$
If $x=\sum_{n=1}^\infty a_nx_n$ then, while computing the norms of $P_n$, we may regard the norm of this sum less than (or equal
to) 1. The inequality (6) gives thus that $\|a_kx_k\|_{A_1}\le 2C$ for any $k$. But, as was proposed from the very beginning,
$\|x_k\|_{A_1}=1$ for all $k$, therefore $|a_k|\le2C$ also for any $k$. Returning to the operators $P_n$, we obtain by (4) that,
for every $x$ with $\|x\|_{A_1}\le1$, one has
$$ \|P_nx\|_{B_1}=\|\sum_{k=0}^\infty a_{n+k}x_{n+k}\|_{B_1}\le2C\sum_{k=0}^\infty\|x_{n+k}\|_{B_1}\le2C\|x_n\|_{B_1}
\sum_{k=0}^\infty2^{-k}. \tag 8 $$
Consequently, $\|P_n\|_{A_1\to B_1}$ tends to zero as $n\to\infty$ due to properties of the sequence $\|x_n\|_{B_1}$.

The estimation of $\|P_n\|_{A_0\to B_0}$ can be done, using inequality (7) with $A_1$ replaced by $A_0$. We get
$$ \|P_n\|_{A_0\to B_0}=\sup_{\|x\|_{A_0}\le1}\|\sum_{k=n}^\infty a_kx_k\|_{B_0}\le(C+1)\sup_{\|x\|_{A_0}\le1}\|\sum_{k=1}^
\infty a_kx_k\|_{A_0}  \tag 9 $$
(the inequality appears due to embedding $A_0\hook B_0$). The last sum in (9) is exactly the norm of $x$ in $A_0$, therefore
$\|P_n\|_{A_0\to B_0}\le C+1$ as required for interpolation (recall that the inequalities (5)--(7) are independent of the
coefficients $a_k$ and we may consider these numbers in (8) and (9) as coordinates of $x$ in different spaces $A_1$ and $A_0$).

For the end of proof, it is enough to use the condition (1) with $T=P_n$, which gives that the norms $\|P_n\|_{A\to B}$
tend to 0 like $\|P_n\|_{A_1\to B_1}$. Recall that the sequence $(x_n)$ is bounded in $A$ by definition. Set $\sup\!\|x_n\|_A=K$. Then
$$ \|x_n\|_B=\|P_nx_n\|_B\le\|P_n\|_{A\to B}\|x_n\|_A\le K\|P_n\|_{A\to B}\rightarrow0\quad\text{as}\ n\to\infty. $$
Consequently, we get that the sequence $(\|x_n\|_B)$ cannot be equivalent to the sequence $(\|x_n\|_A)$.
This contradiction proves the Main Theorem. \qquad \qed \enddemo

{\it Remark.}\ The same proof is applicable to proving interpolation of the ``two-sided" compactness of linear operators
without inequality (1). It is enough to show that $\|P_n\|_{A_0\to B_0}\rightarrow 0$\ like the case of spaces $A_1, B_1$
and to use the inequality $\|T\|_{A\to B}\lesssim\max\{\|T\|_{A_0\to B_0}, \|T\|_{A_1\to B_1}\}$ from the definition of
any interpolation functor. \bs

{\bf Appendix}\ (in accordance with Proposition 1 of [1]). \ Let $T$ be an arbitrary linear operator, acting continuously
from $A_0$ to $B_0$ and compactly from $A_1$ to $B_1$. Define two spaces $Y_i\subset B_i\ (i=0,1)$ consisting of
all $y=Tx,\ x\in A_i$, with norm $\|y\|_{Y_i}=\inf\bigl\{\|x\|_{A_i}\ :\ y=Tx\bigr\}$. As known, every such a space is
called {\it image} of $T$ in the corresponding space $Y_i$ and is a Banach space. Present now $T$
as a composition $T=JS$, where the operator $S$ is defined by $Sx=Tx$ for all $x\in A_0+A_1$ but considered as a bounded
operator from $A_i$ to $Y_i,\ i=0,1$. Consequently $J$ will be the embedding from $Y_i$ into $B_i$ for the same $i=0,1$.
Applying the functor $\F$, we obtain that $S$ acts boundedly from $\F(A_0,A_1)$ into $\F(Y_0,Y_1)$ by standard interpolation
theorems. At the same time, $J$ appears to be compact as an operator from $Y_1$ into $B_1$ due to compactness of the
operator $T$ and thus it is compact as embedding from $\F(Y_0.Y_1)$ into $\F(B_0,B_1)$ in virtue of the Main Theorem, proved
above. In result, $T$ is also compact as an operator from $A$ to $B$. \bs

\enddocument